\documentclass[a4paper,12pt,reqno]{amsart}
\usepackage{amssymb}
\usepackage{amsmath}

\def\la{\langle}
\def\ra{\rangle}
\def\.{\cdot}

\def\n{\nabla}

\def\beq{\begin{equation}}
\def\eeq{\end{equation}}
\def\bea{\begin{eqnarray*}}
\def\eea{\end{eqnarray*}}

\def\f{\varphi}

\def\o{\omega}

\def\L{\Lambda}

\def\r{\end{proof}}

\def \RM{\mathbb{R}}

\def \CM{\mathbb{C}}

\def\d{{\delta}}
\def\dc{d^c}
\def\dec{\delta^c}

\def\es{\,\lrcorner\;}
\def\we{\wedge}
\def\del{{\partial}}
\def\Ric{\mathrm{Ric}}
\def\grad{\mathrm{grad}}
\def\id{\mathrm{id}}
\def\be{\begin{equation}}
\newcommand{\firstline}[2]{\makebox[#1][l]{$\displaystyle{#2}$}&&}
\def\ee{\end{equation}}

\def\pr{\rm{pr}}

\def\Hol{\mathrm{Hol}}
\def\hol{\mathfrak{hol}}
\def\so{\mathfrak{so}}

\def\S{\mathrm{Sym}}

\newtheorem{ede}{Definition}[section]

\newtheorem{epr}[ede]{Proposition}

\newtheorem{ath}[ede]{Theorem}

\newtheorem{elem}[ede]{Lemma}

\newtheorem{ecor}[ede]{Corollary}
 
\title{Twistor forms on K{\"a}hler  manifolds}
\author{Andrei Moroianu and Uwe Semmelmann}
\thanks{The authors are  members of the {\sl European Differential Geometry 
Endeavour} (EDGE), Research Training Network HPRN-CT-2000-00101, 
supported by The European Human Potential Programme;
both authors would like to thank the CNRS for support and the Centre de
Math{\'e}matiques de l'Ecole Polytechnique for hospitality during the
preparation of this work.}
\address{Andrei Moroianu \\ CMAT\\ {\'E}cole Polytechnique \\ UMR 7640 du CNRS
\\ 91128 Palaiseau \\ France}
\email{am@math.polytechnique.fr}

\address{Uwe Semmelmann\\ Mathematisches Institut \\ Universit{\"a}t
M{\"u}nchen  \\ Theresienstr. 39 \\ D-80333 M{\"u}nchen, Germany}
\email{Uwe.Semmelmann@mathematik.uni-muenchen.de}

\begin{document}

\begin{abstract}
Twistor forms are a natural generalization of conformal 
vector fields on Riemannian manifolds. They are defined as 
sections in the kernel
of a conformally invariant first order differential operator. We study 
twistor forms on compact K{\"a}hler manifolds and give a complete 
description up to special forms in the middle dimension. 
In particular, we show that they are closely related to Hamiltonian 
2-forms. This provides the first examples of compact K{\"a}hler
manifolds with non--parallel twistor forms in any even degree.
\vspace{0.1cm}

\end{abstract}

\maketitle

\section{Introduction}

{\it Killing vector fields} are important objects in Riemannian
geometry. They are by definition infinitesimal isometries, i.e.~their
flow preserves a given metric. The existence of 
Killing vector fields determines the degree of symmetry of the manifold.  
Slightly more generally one can consider
{\it conformal vector fields}, i.e. vector fields whose flows
preserve a conformal class of metrics. The covariant derivative
of a vector field can be seen as a section of the tensor product
$\L^1 M\otimes TM$ which is isomorphic to $\L^1M\otimes \L^1M$. This
tensor product decomposes under the action of $O(n)$ as
$$\L^1M\otimes \L^1M\cong \RM\oplus\L^2M\oplus S^2_0M,$$
where $S^2_0M$ is the space of trace--free symmetric 2--tensors,
identified with the {\it Cartan product} of the two copies of
$\L^1M$. A vector field $X$ is a
conformal vector field if and only if the projection on $S^2_0M$ of
$\n X$ vanishes.

More generally, the tensor product $\L^1M\otimes \L^pM$ decomposes
under the action of $O(n)$ as
$$\L^1M\otimes \L^pM\cong \L^{p-1}M\oplus\L^{p+1}M\oplus \L^{p,1}M,$$
where again $\L^{p,1}M$ denotes the Cartan product. As natural
generalizations of conformal vector fields, {\it twistor forms} are
defined to be $p$--forms $\psi$ such that the projection of $\n\psi$
onto $\L^{p,1}M$ vanishes.

Coclosed twistor $p$--forms are called {\it Killing forms}. For $p=1$ 
they are dual to Killing vector fields. Note that parallel forms
are trivial examples of twistor forms.

Killing forms, as generalization of Killing vector fields, 
were introduced by K.~Yano in \cite{yano}. Twistor forms were
introduced later on by
S.~Tachibana \cite{ta3}, for the case of 2--forms,
and by T.~Kashiwada \cite{ka}, \cite{ta4} in the general case.

The projection $\L^1M\otimes \L^pM\to \L^{p,1}M$  defines a first order 
differential operator $T$, which was already studied in the 
context of Stein--Weiss operators (c.f. \cite{branson}). As 
forms in the kernel of $T$, twistor forms  are very similar to twistor
spinors in spin geometry. We will give an explicit construction
relating these two objects in Section 2.

The special interest for twistor forms in the physics literature
stems from the fact that they can be used to define quadratic first 
integrals of the geodesic equation, i.e. functions which are constant 
along geodesics. Hence, they can be used to integrate the equation
of motion, which was done for the first time by R. Penrose and
M. Walker in \cite{penrose}.
More recently Killing forms and twistor 
forms have been successfully applied to define symmetries of field 
equations (c.f.~\cite{be1}, \cite{be2}).

Despite this longstanding interest in Killing forms there are
only very few global results on the existence or non--existence
of twistor forms on Riemannian manifolds. In this article we will
give a description of twistor forms on compact K{\"a}hler
manifolds and show how they are related to the Hamiltonian 2--forms
which were recently introduced in \cite{acg1}.

The first result for twistor forms on compact K{\"a}hler manifolds
was obtained in \cite{yama3}, where it is proved that any Killing
form has to be parallel. Some years later S.~Yamaguchi et al. stated
in \cite{yama} that, with a few exceptions for small degrees and in 
small dimensions, any twistor form  on a compact K{\"a}hler 
manifold has to be parallel. Nevertheless, it turns out that their
proof contains several serious gaps.

In fact, we will show that there are examples of compact K{\"a}hler
manifolds having non--parallel twistor forms in any even degree.
The simplest examples can be found on the complex projective space.
More precisely we will show that on a compact K{\"a}hler manifold $M^{2m}$
any twistor $p$--form with $p\ne m$ can (modulo parallel forms) be written as 
$L^{k-1}\varphi +  L^kf$,
where $\varphi\in\L^2 M$ and $f\in C^\infty M$ are, up to some
constants depending on $m$ and $k$, the primitive part and the trace
of a {\it Hamiltonian} 2--form (here $\,L\,$ denotes wedging with the K{\"a}hler 
form). Conversely, any Hamiltonian 2--form
can be used to construct twistor forms. 
Hamiltonian 2--forms were recently studied in \cite{acg1} and 
\cite{acg}. They appear in many different areas of K{\"a}hler
geometry. In particular they arise on weakly Bochner--flat K{\"a}hler manifolds
and on K{\"a}hler manifolds which are conformally Einstein.

The case of twistor 
$m$--forms is more complicated. They are either of the form
$L^{k-1}\varphi$, where $\varphi$ is the primitive part of a
Hamiltonian 2--form,
or they belong to a special class of forms, whose existence is still
open for the time being.

\newpage


\section{Twistor forms on Riemannian manifolds}

Let $( V,\,\la\cdot,\cdot\ra)$ be an $n$--dimensional Euclidean vector space. 
The tensor product  $V^*\otimes\Lambda^pV^*$ has the following
$O(n)$--invariant decomposition:
\bea\label{deco1}
V^*\otimes\Lambda^pV^*
&\cong&
\Lambda^{p-1}V^* 
\oplus \Lambda^{p+1}V^* 
\oplus 
\Lambda^{p,1}V^* 
\eea
where $\Lambda^{p,1}V^*$ is the intersection of the kernels of wedge
and inner product maps, which can be identified with the Cartan
product of $V^*$ and $\Lambda^pV^*$. This decomposition  immediately
translates to Riemannian manifolds $\,(M^n,\,g)$, where we have 
\begin{equation}\label{deco}
T^*M\otimes\Lambda^pT^*M
\;\cong\;
\Lambda^{p-1}T^*M
\oplus
\Lambda^{p+1}T^*M
\oplus
\Lambda^{p,1}T^*M
\end{equation}
with $\Lambda^{p,1}T^*M$ denoting the vector bundle corresponding to 
the vector space 
 $\Lambda^{p,1}$. The covariant derivative $\nabla \psi$ 
of a $p$--form $\psi$ 
is a section of $\;T^*M\otimes\Lambda^pT^*M$. Its projections onto
the summands $\,\Lambda^{p+1}T^*M\,$ and $\,\Lambda^{p-1}T^*M\,$
are just the differential
$d\psi$ and the codifferential $\d \psi$.
Its projection onto the third summand $\,\Lambda^{p,1}T^*M\,$ defines a
natural first order differential operator $T$, called the 
{\it twistor operator}. 
The twistor operator
$ 
T:\Gamma(\Lambda^p T^*M)\,\rightarrow \,\Gamma(\Lambda^{p,1}T^*M) 
\subset 
\Gamma(T^*M\otimes\Lambda^pT^*M)
$
is given for any vector field $X$ by the following formula
$$
[\,T\psi\,]\,(X)
\;:=\;
[\pr_{\Lambda^{p,1}}(\nabla \psi)]\,( X)
\;=\;
\nabla_X\, \psi
\;-\;
\tfrac{1}{  p+1}\, X \es d\psi
\;+\;
\tfrac{1}{ n-p+1}\, X\,\wedge\,\d \psi \ .
$$
Note that here, and in the remaining part
of this article, we identify vectors and 1--forms using the metric.

The twistor operator $T$ is a typical example of a so-called
Stein--Weiss operator and it was in this context already considered 
by T.~Branson in \cite{branson}. It particular it was shown
that $T^*T$ is elliptic, which easily follows from computing 
the principal symbol.
Its definition is also similar to the definition of the twistor
operator in spin geometry. The
tensor product between the spinor bundle and the cotangent bundle
decomposes under the action of the spinor group into the sum
of the spinor bundle and the kernel of the Clifford
multiplication. The (spinorial) 
twistor operator is then defined as the projection of the covariant derivative
of a spinor onto the kernel of the Clifford multiplication.

\begin{ede}
A $p$--form $\,\psi\,$ is called a {\it twistor $p$--form} if and only if 
$\,\psi\,$ is in the kernel of $\,T$,~i.e. 
if and only if $\,\psi\,$ satisfies 
\begin{equation}\label{killing}
\nabla_X\,\psi\;=\;
\tfrac{1}{  p+1}\,X\es d\psi \;-\;
\tfrac{1}{  n-p+1}\, X\,\wedge\,\d \psi\ ,
\end{equation}
for all vector fields $X$. 
If the $p$--form $\psi$ is in addition coclosed, it is called a {\it Killing $p$--form}. 
This is equivalent to $\nabla\psi\in\Gamma(\Lambda^{p+1}T^*M)$ or to 
$X \es \nabla_X \psi = 0$ for any vector field $X$.
\end{ede}
It follows directly from the definition that the Hodge 
star-operator $\ast$ maps twistor $p$--forms into twistor
$(n-p)$--forms. In particular, it interchanges closed and coclosed
twistor forms.

In the physics literature, equation~(\ref{killing}) defining a twistor 
form is often called the {\it Killing--Yano equation}. The terminology
{\it conformal Killing forms} to denote twistor forms is also used. 
Our motivation for using the name {\it twistor form} is not only
the similarity of its definition to that of twistor spinors in
spin geometry, but also the existence of a direct relation between
these objects. We recall that a {\it twistor spinor} on a Riemannian spin
manifold is a section $\,\psi\,$ of the spinor bundle lying in the kernel
of the (spinorial) twistor operator. Equivalently, $\psi$ satisfies for
all vector fields $X$ the equation
$
\nabla_X \psi \,=\, -\,\tfrac{1}{n}\,X \cdot D\psi
$,
where $D$ denotes the Dirac operator. Given two such twistor spinors,
$\psi_1$ and $\psi_2$, we can introduce $k$--forms $\omega_k$, which
on tangent vectors $X_1, \ldots , X_k$ are defined as
$$
\omega_k(X_1, \ldots, X_k)
:=
\la (X_1 \wedge \ldots \wedge X_k)\cdot \psi_1, \psi_2 \ra \ .
$$
It is well--known that for $k=1$ the form $\omega_1$ is dual to
a conformal vector field. Moreover, if $\psi_1$ and $\psi_2$ are
Killing spinors the form $\omega_1$ is dual to a Killing
vector field. More generally we have

\begin{epr} {\em (cf. \cite{uwe})}
Let $(M^n, g)$ be a Riemannian spin manifold with twistor 
spinors $\psi_1$ and $\psi_2$. Then for any $k$ the associated
$k$--form $\omega_k$ is a twistor form. 
\end{epr}

We will now give an important integrability condition which
characterizes twistor forms on compact manifolds. A similar 
characterization was obtained in~\cite{ka}. We first 
obtain two Weitzenb{\"o}ck formulas by differentiating 
the equation defining the twistor operator.
\begin{eqnarray}
\nabla^*\nabla\psi 
&=& \label{weiz1}
\tfrac{1}{ p+1}\,
\d d\, \psi\;\;+\;\; \tfrac{1}{  n-p+1}\,d \d \, \psi\;\;+\;\;T^*T\,
\psi \ ,
\\[1.5ex]
q(R)\,  \psi 
&=&\label{weiz2}
\tfrac{p}{ p+1}\,\d d\, \psi
\;\;+\;\;\tfrac{n-p}{ n-p+1}\,d \d \, \psi\;\;-\;\;T^*T\, \psi \ ,
\end{eqnarray}
where $q(R)$ is the curvature expression appearing in the classical
Weitzenb{\"o}ck formula for the Laplacian on $p$--forms:
$\Delta  = \d d\,+\,d\d  = \nabla^*\nabla \,+\, q(R) $. It is  the symmetric
endomorphism of the bundle of differential forms defined by
\begin{equation}\label{qr}
q(R)\;=\;\sum\,e_j\,\wedge\,e_i \es R_{e_i,e_j},
\end{equation}
where $\,\{e_i\}\,$ is any local orthonormal frame and $R_{e_i,e_j}$
denotes the curvature of the form bundle. On forms of degree one and
two one has an explicit expression for the action of
$q(R)$,~e.g. if $\xi$ is any 1--form, then $q(R)\,\xi = \Ric(\xi)$.
In fact it is possible to define $q(R)$ in a more general
context. For this we first rewrite equation~(\ref{qr}) as
$$
q(R)
\;=\;
\sum_{i < j}\,
(e_j\,\wedge \,e_i \es \,-\; e_i\,\wedge \,e_j\es)\, R_{e_i,e_j}
\;=\;
\sum_{i < j}\,
(e_i \we e_j)\bullet R(e_i \we e_j )\bullet
$$
where the Riemannian curvature $R$ is considered as element of
$\, \S^2(\L^2 T_pM)\,$ and $\,\bullet\,$ denotes the standard representation
of the Lie algebra $\,\so(T_pM)\cong\L^2 T_pM \,$ on the space of
$p$--forms. Note that we can replace $\,e_i \we e_j\,$ by any orthonormal
basis  of $\so(T_pM)$.
Let $(M, g)$ be a Riemannian manifold with holonomy group $\Hol$. Then
the curvature tensor takes values in the Lie algebra 
$\,\hol \, $ of the 
holonomy group,~i.e. we can write $q(R)$ as 
$$
q(R)
\;=\;
\sum \, \omega_i \bullet R(\omega_i)\bullet \qquad \in \;\S^2(\hol)
$$
where $\{\omega_i\}$ is any orthonormal basis of $\, \hol\,$ and $\,\bullet\,$
denotes form representation restricted to the holonomy group.
Writing the bundle endomorphism $q(R)$ in this way has two
immediate consequences: we see that $q(R)$ preserves any
parallel subbundle of the form bundle and it is clear that
by the same definition $q(R)$ gives rise to a symmetric endomorphism
on any associated vector bundle defined via a representation of
the holonomy group. 

Integrating the second Weitzenb{\"o}ck formula (\ref{weiz2}) yields
the following integrability condition for twistor forms. 
\begin{epr}\label{integrabl}
Let $(M^n,\,g)$ a compact Riemannian manifold. Then a $p$--form $\psi$ 
is a twistor $p$--form, if and only if
\begin{equation}\label{char}
q(R)\,  \psi 
\;\;=\;\;
\tfrac{p}{ p+1}\,\d d\, \psi\;\;+\;\;\tfrac{n-p}{ n-p+1}\,d \d \, \psi \ .
\end{equation}
\end{epr}

For coclosed forms, Proposition \ref{integrabl} is a generalization of 
the well--known characterization of Killing vector fields on compact
manifolds, as divergence free vector fields in the kernel of 
$\,\Delta - 2\,\Ric$. In the general case, it can be 
reformulated as
\begin{ecor}\label{integrab}
Let $(M^n,\,g)$ a compact Riemannian manifold with a coclosed $p$--form
$\psi$. Then $\psi$ is a Killing form if and only if
$$
\Delta \, \psi
\;\;=\;\;
\frac{p+1}{p}\,q(R)\,  \psi \ .
$$
\end{ecor}
One has similar characterizations for closed twistor forms and for
twistor m--forms on 2m--dimensional manifolds.

An important property of the equation defining 
twistor forms is its conformal invariance (c.f.~\cite{be2}).
We note that the same is true for the twistor equation in spin
geometry. The precise formulation for twistor forms is the following.
Let  $\psi$ be a twistor form on a Riemannian manifold $(M, g)$. Then 
${\widehat \psi}:= e^{(p+1)\lambda}\psi$ is a twistor $p$--form with 
respect to the conformally equivalent metric ${\hat g}:= e^{2\lambda}g$.
Parallel forms are obviously in the kernel of the twistor operator,
hence they are twistor forms.  Using the conformal invariance we see that 
with any parallel form $\,\psi$, also the form
$\,{\widehat \psi}:= e^{(p+1)\lambda}\,\psi\,$ is a twistor
$\,p$--form with respect to the conformally equivalent
metric $\,{\widehat g}:= e^{2\lambda} g$. The form $\,{\widehat
\psi}\,$ is in general not parallel.

The first non trivial examples of twistor forms were found on the
standard sphere. Here it is easy to show that any twistor $p$--form
is a linear combination of eigenforms for the minimal eigenvalue
of the Laplace operator on closed resp. coclosed $p$--forms. It  
is shown in \cite{uwe} that the number of linearly independent
twistor forms on a connected Riemannian manifold is bounded from
above by the corresponding number on the standard sphere.
Other examples exist on Sasakian, nearly K{\"a}hler and on
weak $G_2$--manifolds. In some sense they are all related to the
existence of Killing spinors. Up to now these are more or less all known
global examples of twistor forms. 

In section~\ref{kaehler} we will see that there are many new examples 
on compact K{\"a}hler manifolds. The examples include the complex
projective space and  Hirzebruch surfaces.


\section{Twistor forms on K{\"a}hler manifolds}\label{kaehler}

In this section we will consider twistor forms on a compact
K{\"a}hler manifold $\,(M,\,g,\,J)\,$ of dimension $n=2m$. The case
of forms in the middle dimension,~i.e. forms of degree $m$, is 
somewhat special and will be treated in the next section.

On a complex manifold the differential splits as $d = \del +{\bar \del}$.
Moreover, one has the following real differential operator
$$
d^c 
\;=\;
i\,({\bar \del} \,-\, \del)
\;=\;
\sum\, Je_i \we \nabla_{e_i} \ .
$$  
The formal adjoints of $d$ and $\dc$ are denoted by $\d$ and $\d^c$.  If $\,\omega\,$ is the K{\"a}hler form   
then  $\Lambda$ denotes the  contraction with $\omega$ and $L$ the
wedge product with $\omega$. Another important operator acting on
forms is the natural extension of the complex structure $J$ on $-p$--forms, defined
as follows~:
$$Ju:=\sum J(e_i)\wedge e_i \es \, u\qquad\forall u\in\L^pM.$$
These operators satisfy the
following fundamental commutator relations  
\bea
d^c 
&=& -\,[\,\d ,\,L\,]
\;=\;
-\,\,[\,d ,\,J \,] \, ,
\qquad
\d^c
\;=\;
\,[\, d,\, \Lambda\,]
\;=\;
-\,[\,\d ,\,J \,] \ ,
\\[1.5ex]
d 
&=&
 \,[\, \d^c,\, L\,]
\;=\;
 \,[\, d^c,\,J \,] \ ,
\qquad
\d
\;=\;
- \,[\,d^c ,\, \Lambda\,]
\;=\;
 \,[\,\d^c ,\,J \,] \ .
\eea
In addition the following commutators resp. anti--commutators vanish
\bea
0
&=&
\,[\,d ,\,L\,]
\;=\;
\,[\,d^c ,\,L \,]
\;=\;
\,[\,\d ,\, \Lambda\,]
\;=\;
\,[\, \d^c,\, \Lambda\,]
\;=\;
\,[\,\Lambda ,\,J \,]
\;=\;
\,[\,J ,\, \ast\,] \ ,
\\[1.5ex]
0
&=&
\d d^c \;+\; d^c \d 
\;=\;
dd^c + d^cd
\;=\;
\d\d^c \;+\; \d^c\d
\;=\;
d\d^c \;+\;\d^cd \ .
\eea
Using these relations we are now able to derive several consequences
of the twistor equation on K{\"a}hler manifolds. We start by computing 
$\d^c u$ for a twistor $p$--form $u$.
\begin{eqnarray}
\d^c u
&=&
-\,\sum\,
Je_i \es \nabla_{e_i} \,u 
\;=\;
-\,\sum\,
Je_i \es
(
\tfrac{1}{p+1}\,e_i \es du \;-\; \tfrac{1}{n-p+1} \, e_i \we \d u
)
\nonumber
\\[1.5ex]
&=&
-\,\tfrac{2}{p+1}\,\Lambda(du)
\;+\;
 \tfrac{1}{n-p+1} \,J(\d u)
\nonumber
\\[1.5ex]
&=&
-\,\tfrac{2}{p+1}\,d\Lambda(u)
\;+\;
\tfrac{2}{p+1}\,\d^c\,u
\;+\;
 \tfrac{1}{n-p+1} \,J(\d u)
\nonumber
\\[1.5ex]
&=&
-\,\tfrac{2}{p-1}\,d\Lambda(u)
\;+\; 
\tfrac{p+1}{(p-1)(n-p+1)} \,J(\d u) \ . \label{zwei}
\end{eqnarray}
From this we immediately obtain
\begin{eqnarray}
\d^c\,d\,u &=& -\,d\,\d^c \,u
\;=\;
-\,\tfrac{p+1}{(p-1)(n-p+1)} 
\left(
J(d\d u)
\;+\;
d^c\d \,u
\right) \ , \label{drei}
\\[1.5ex]
\d\d^c \, u 
&=&
-\,\d^c\d u 
\;=\;
-\,\tfrac{2(n-p+1)}{(n-p)(p+1)}\,\Lambda(\d du) \label{eins} \ .
\end{eqnarray}
A similar calculation for the operator $\,d^c\,$ leads to
\begin{eqnarray}
d^c\,u
&=&
-\,\sum\, e_i \we \nabla_{Je_i}\,u
\;=\;
-\,\sum\, e_i \we
(
\tfrac{1}{p+1}\,Je_i \es du
\;-\;
\tfrac{1}{n-p+1}\,Je_i \we \d u
)
\nonumber
\\[1.5ex]
&=&
\tfrac{1}{p+1}\,J(du)
\;+\;
\tfrac{2}{n-p+1}\,L(\d u)
\;=\;
\tfrac{1}{p+1}\,J(du)
\;+\;
\tfrac{2}{n-p+1}\,\d (L u)
\;+\;
\tfrac{2}{n-p+1}\,d^cu
\nonumber
\\[1.5ex]
&=&
\tfrac{n-p+1}{(p+1)(n-p-1)}\,J(du)
\;+\;
\tfrac{2}{n-p+1}\,\d (L u) \ . \label{null}
\end{eqnarray}
We apply this to conclude
\begin{eqnarray}
\d d^c\,u
&=&
\tfrac{n-p+1}{(p+1)(n-p-1)}\,
(\,J \d d u \;-\; \d^cdu
\,) \ . \label{vier}
\end{eqnarray}

As a first step we will prove that for a twistor $p$--form $\,u\,$
the form $\,\Lambda J(u)\,$ has to be parallel. Here we can
assume $\,p \neq 2$, since otherwise there is nothing to 
show. It will be convenient to introduce for a moment the 
following notation
\bea
x 
&:=& 
J \Lambda (d\d  u), 
\qquad
a 
\;:=\;
d\d (J \Lambda u), 
\qquad 
\alpha 
\;:=\;
\d^cd(\Lambda u),
\\[1.5ex]
\qquad  
y 
&:=&
J \Lambda (\d d u) ,
\qquad
b 
\;:=\;
\d d (J \Lambda u),
\qquad 
\beta 
\;:=\; \d d^c(\Lambda u) \ .
\eea

Using the various commutator rules we will derive several equations
relating the quantities $\, x,\,y,\,a,\,b,\,\alpha\,$ and $\,\beta$. 
We start with
\bea
x
&=&
J \Lambda ( d\d u)
\;=\;
J( d\d  \Lambda u \;-\; \d^c \d  u)
\\[1.5ex]
&=&
d\d (J\Lambda u) \;+\; d\d^c (\Lambda u) \;+\; d^c\d (\Lambda u)
\;-\; J\d^c \d u 
\\
[1.5ex]
&=&
a \;-\; \alpha \;-\; \beta \,-\,
\tfrac{2(n-p+1)}{(n-p)(p+1)}\,y \ ,
\eea
where we also used equation~(\ref{eins}) to replace the summand
$\,J\d^c \d u$. Similarly we have
\bea
y
&=&
J\Lambda(\d d u)
\;=\;
J(\d d(\Lambda u) \;-\; \d \d^c u)
\\[1.5ex]
&=&
\d d (J\Lambda u) \;+\; \d d^c(\Lambda u) \;+\; \d^cd(\Lambda u)
 \;-\; J\d\d^cu 
\\[1.5ex]
&=&
b \;+\; \beta \;+\; \alpha \,+\,
\tfrac{2(n-p+1)}{(n-p)(p+1)}\,y \ .
\eea
Adding the two equations for $x$ resp. $y$, we obtain that
$\,x+y=a+b$,~i.e. we can express $x$ and $y$
in terms of $a,\,b,\,\alpha$ and $\,\beta\,$ as
\begin{eqnarray}
y
&=&
\epsilon_1\,(b\;+\;\alpha\;+\;\beta),
\qquad \qquad 
\epsilon_1 = \tfrac{(n-p)(p+1)}{(n-p)(p+1)-2(n-p+1)}\ ,
\label{y}
\\[1.5ex]
x
&=&
a\;+\;b\;-y\;
\;=\;
a\;+\;(1\,-\,\epsilon_1)\,b\;-\;\epsilon_1\,(\alpha\,+\beta)
\label{x} \ .
\end{eqnarray}
Note that $\,\epsilon_1\,$ is well defined since the denominator
may vanish only in the case $\,(n,p)=(4,2)$, which we already excluded.
Next, we contract equation~(\ref{drei}) with the K{\"a}hler form
to obtain an equation for $\,\alpha$.
\bea
\alpha
&=&
\d^c d (\Lambda u) \;=\; \Lambda (\d^c d u)
\\[1.5ex]
&=&
-\,\tfrac{p+1}{(p-1)(n-p+1)} 
\left(
\Lambda J(d\d u)
\;+\;
\Lambda (d^c\d\,u)
\right)
\\[1.5ex]
&=&
-\,\tfrac{p+1}{(p-1)(n-p+1)}
\,(x\;-\;\beta) \ .
\eea
Contracting equation~(\ref{vier}) with the K{\"a}hler form we obtain
a similar expression for $\,\beta$.
\bea
\beta 
&=&
\d d^c(\Lambda u) \;=\; \Lambda(\d d^c u)
\\[1.5ex]
&=&
\tfrac{n-p+1}{(p+1)(n-p-1)}\,
(\,\Lambda J \d d u \;-\; \Lambda(\d^cdu)
\,) 
\\[1.5ex]
&=&
\tfrac{n-p+1}{(p+1)(n-p-1)}\,(y \;-\; \alpha) \ .
\eea
Replacing $x$ and $y$ in the equations for $\alpha$ resp. $\beta$
leads to
\begin{eqnarray*}
a\;+\;(1\,-\,\epsilon_1)\,b\;-\;\epsilon_1\,(\alpha\,+\beta)
&=&
\epsilon_2\,\beta \;+\; \alpha,
\qquad \qquad
\epsilon_2 = \tfrac{(p-1)(n-p+1)}{p+1}
\\[1.5ex]
\epsilon_1\,(b\;+\;\alpha\;+\;\beta)
&=&
\epsilon_3\,\beta \;+\; \alpha
\qquad \qquad
\epsilon_3 = \tfrac{(p+1)(n-p-1)}{n-p+1} \ .
\end{eqnarray*}
Finally we obtain two equations only involving $\,a,\,b,\,\alpha\,$
and $\,\beta$,
\begin{eqnarray}
a\;+\;(1\,-\,\epsilon_1)\,b
&=&
(\epsilon_1\,-\,\epsilon_2)\,\alpha \;+\; (\epsilon_1\,-\,1)\,\beta
\label{stern}
\\[1.5ex]
\epsilon_1\,b
&=&
(1\,-\,\epsilon_1)\,\alpha \;+\; (\epsilon_3\,-\,\epsilon_1)\,\beta
\label{zweistern} \ .
\end{eqnarray}
Considering equation~(\ref{zweistern}) we note that $b$ and $\beta$
are in the image of $\,\d$, whereas $\,\alpha\,$ is in the image
of $\,d$. Hence, since the coefficient $\,(1-\epsilon_1)\,$ is
obviously different from one, we can conclude that $\,\alpha\,$ has
to vanish. Moreover, we have the equation 
$\,\epsilon_1 b = (\epsilon_3 - \epsilon_1)\beta$. Applying the
same argument in equation~(\ref{stern}) we obtain that $\,a\,$
has to vanish and that $\,(1-\epsilon_1)b = (\epsilon_1-1)\beta$.
Combining the two equations for $\,b\,$ and $\,\beta\,$ we see
that also these expressions have to vanish. Taking the scalar product
of $\,a\,$ and $\,b\,$ with $\,J\Lambda u\,$ and integrating over $M$
yields that $\,J\Lambda u\,$ has to be closed and coclosed, hence
harmonic. 

\begin{epr}\label{jlu}
Let $(M,\,g,\,J)$ be a compact K{\"a}hler manifold. Then for any
twistor form $\,u\,$ the form $\,\Lambda J u \,$ is parallel.
\end{epr}
\proof
Let $\,u\,$ be a twistor $p$--form. We have already shown that
$\,J\Lambda u\,$ is harmonic, i.e. $\,\Delta( J\Lambda u) =0$.
Since $\,\Delta = \nabla^*\nabla + q(R)$, the proof of the
proposition would follow from $q(R)(J\Lambda u)=0$. For a 
twistor form $\,u\,$ we have the Weitzenb{\"o}ck formula~(\ref{weiz2}) and
the argument in Section 2 shows that the curvature endomorphism $\,q(R)\,$
commutes with $\Lambda$ and $J$. Hence, we obtain 
$$
q(R) (J \Lambda u)
\;=\;
\tfrac{p}{p+1}\,x \;+\; \tfrac{n-p}{n-p+1}\,x \ ,
$$
with the notation above. But with 
$\,a,\,b,\,\alpha\,$ and $\,\beta\,$ also $\,x\,$
and $\,y\,$ have to vanish. Hence, $\,q(R)(J\Lambda u)=0\,$
and the Weitzenb{\"o}ck formula for $\,\Delta \,$ and integration 
over $M$ prove that $\,J \Lambda u\,$ has to be parallel.
\qed

In the case where the degree of twistor form is different from
the complex dimension, this proposition has an important corollary.

\begin{ecor}\label{jin}
Let $(M,\,g,\,J)$ be a compact K{\"a}hler manifold of complex
dimension $m$. Then for any twistor $p$--form $\,u$, with $\,p\neq m$,
$Ju$ is parallel.
\end{ecor}
\proof
We first note that the Hodge star operator $\ast$ commutes with
the complex structure $J$ and with the covariant derivative. Moreover, it
interchanges $L$ and $\Lambda$. If $\,u\,$ is a twistor $p$--form, then
$\,\ast u\,$ is a twistor form, too, and by Proposition \ref{jlu}, 
$\,\Lambda J \ast u\,$ is parallel. Hence, also 
$\,\ast \Lambda J \ast u\,$ is parallel. But 
$\,\ast \Lambda J \ast u = \pm L J u\,$ and it follows that
$\,L J u \,$ is parallel. Since the contraction with the 
K{\"a}hler form commutes with the covariant derivative we conclude
that $\,\Lambda L J u \,$ is parallel as well. Thus $(m-p)Ju=\L L
Ju-L\L Ju$ is parallel.
\qed

In the remaining part of this section we will investigate twistor
forms $\,u\,\,$ for which $\,Ju\,$  
is parallel. The results below are thus valid for all twistor $p$--forms
with $\, p\neq m$, but also for $J$--invariant twistor $m$--forms, a
fact used in the next section.

\begin{epr}\label{differentials}
Let $(M,\,g,\,J)$ be a compact K{\"a}hler manifold of dimension $n=2m$. 
Then any twistor $p$--form $\,u\,$ for which
$\,Ju\,$ is parallel satisfies the equations
\begin{eqnarray*}
\d^c u 
&=&
\mu_1 \, d\Lambda u,
\qquad
d^c u 
\;=\;
\mu_2 \, \d L u,
\qquad
\d u
\;=\;
-\,\mu_1\,d^c \Lambda u, 
\qquad
d u \;=\; -\, \mu_2 \, \d^c L u \ ,
\end{eqnarray*}
with constants
$
\;
\mu_1 := -\,\tfrac{2(n-p+1)}{(p-1)(n-p)-2}
\;
$
and
$
\;
\mu_2 := \tfrac{2(p+1)}{p(n-p-1)-2}
$.
\end{epr}
Note that this proposition is also valid for forms in the middle dimension.
\proof
The formula for $\,\d^cu\,$ follows from equation~(\ref{zwei}) after
interchanging $J$ and $\d$ and using the assumption that $\,Ju\,$
is parallel. Because of $\,J\d = \d J + \d^c\,$ we obtain
$$
\d^c u\,
(1-\tfrac{p+1}{(p-1)(n-p+1)})
\;=\;
-\,\tfrac{2}{p-1}\,d\Lambda(u)
\;+\; 
\tfrac{p+1}{(p-1)(n-p+1)} \,\d (J u)
\;=\;
-\,\tfrac{2}{p-1}\,d\Lambda(u) \ .
$$
Since $\,J \Lambda u\,$ is also parallel, the formula for
$\,\d u\,$ follows by applying $J$ to the above equation. Similarly we
obtain the expression for $\,d^cu\,$ by using equation~(\ref{null})
and the relation $\,Jd=dJ + d^c$. This yields
$$
d^c u\,
(1-\tfrac{n-p+1}{(p+1)(n-p-1)})
\;=\;
\tfrac{n-p+1}{(p+1)(n-p-1)}\,d (J u)
\;+\;
\tfrac{2}{n-p-1}\,\d(L u)
\;=\;
\tfrac{2}{n-p-1}\,\d(L u) \ .
$$
Applying $J$ to this equation yields the formula for $\,d u\,$.
\qed

As a first application of Proposition~\ref{differentials} we will show
that the forms
$\,du\,$  and $\,\d u\,$ are eigenvectors of the operator 
$\,\Lambda L$.

\begin{elem}\label{myh}
Let $\,u\,$ be a twistor $p$--form with $\,Ju\,$ parallel. Then
$\,du\,$ and $\,\d u\,$ satisfy the equations
$$
\Lambda  L (du) 
\;=\;
\tfrac{1}{4}\,(n-p-2)\,(p+2)\, du \ ,
\qquad\quad
\Lambda  L (\d u) 
\;=\;
\tfrac{1}{4}\,(n-p)\,p\,\d u \ .
$$
\end{elem}
\proof
From Proposition~\ref{differentials} we have the equation
$\,du = - \mu_2\, \d^c L u = - \mu_2 L \d^c u - \mu_2 u$. This implies
\bea
(1\,+\,\mu_2)\,du
&=&
-\,\mu_2 L (\d^c u)
\;=\;
-\,\mu_1\,\mu_2\, L d \Lambda u
\;=\;
-\,\mu_1\,\mu_2\, d L \Lambda u
\\[1.5ex]
&=&
-\,\mu_1\,\mu_2\, d ( \Lambda L u \,-\, (m-p)\,du)
\\[1.5ex]
&=&
-\,\mu_1\,\mu_2\,
(\d^c Lu \,+\, \Lambda L u \,-\, (m-p)\,du)
\\[1.5ex]
&=&
\mu_1\, du \;-\; \mu_1\,\mu_2\,\Lambda L (du)
\;+\;
\mu_1\,\mu_2\,(m-p)\,du \ .
\eea
Collecting the terms with $\,du\,$ we obtain
$$
\Lambda L(du)
\;=\;
-\,\tfrac{1\,+\,\mu_2\,-\,\mu_1\,-\,\mu_1\mu_2(m-p)}{\mu_1\mu_2}\,du
$$
where $m$ denotes the complex dimension. Substituting the
values for $\,\mu_1\,$ and $\,\mu_2$ given in 
Proposition~\ref{differentials}, the formula for 
$\,\Lambda L (du)\,$ follows after a straightforward calculation.

\medskip

Similarly we could prove the formula for $\,\Lambda L(\d u)\,$
by starting from the equation $\,d^cu= \mu_2\, \d Lu$. Nevertheless
it is easier to note that $\,\ast u\,$ is again a twistor form
(of degree $n-p$).
Hence, $\,4 \Lambda L(d\ast u)=(p-2)(n-p+2) d\ast u\,$ and
we have
$$
(p-2)(n-p+2) \, \ast d \ast u
\;=\;
4 \ast \Lambda L(d\ast u)
\;=\;
4 \, L \Lambda \ast d \ast u 
\;=\;
4 \, \Lambda L \ast d \ast u \;-\; 4\,(m-p+1) \ast d \ast u \ .
$$
Hence,
$$
4\,\Lambda L (\d u)
\;=\;
[(p-2)(n-p+2) \,+\,4\,(m-p+1) ]\, \d u
\;=\;
(n-p)\,p\,\d u \ .
$$
\qed

In order to use the property that $\,du\,$ resp. $\,\d u\,$
are eigenvectors of $\,\Lambda L\,$ we still need the
following elementary lemma.

\begin{elem}\label{elementary}
Let $ \,(M^{2m},g, J)$ be a K{\"a}hler manifold and let $ \,\alpha \,$ be a 
$ \,p$--form on $ \,M$, then
$$
[\,\Lambda, \, L^s\,]\,\alpha
\;=\;
s\,(m-p-s+1)\,L^{s-1}\alpha \ .
$$
Moreover, if $\,\alpha\,$ is a primitive $p$--form, i.e. if 
$\,\Lambda(\alpha)=0$, then $\,\alpha\,$ satisfies in addition
the equation
$$
\Lambda^r \, L^s \, \alpha
\;=\;
\tfrac{s!(m-p-s+r)!}{(s-r)!(m-p-s)!}\,L^{s-r}\alpha \ ,
$$
where $r, s$ are any integers.
\end{elem}
\proof
In the case $r=1$ we prove the formula for the commutator of $\Lambda$ 
and $L$ by induction with respect to $s$. 
For $s=1$ it is just the well--known commutator relation for
$\Lambda$ and $L$. Assume now that we know the formula for $s-1$, then
\bea
[\,\Lambda, \, L^s\,]\,\alpha
&=&
(\,\Lambda \circ L^s - L^s \circ \Lambda\,)\,\alpha
\;=\;
\Lambda \circ L^{s-1} \,(L \alpha)
\;-\;
L^{s-1}\,(L \circ \Lambda )\, \alpha
\\[1.5ex]
&=&
L^{s-1} \circ \Lambda \,(L \alpha)
\;+\;
(s-1)\,(m-(p+1)-(s-1)+1)\,L^{s-2}\,(L\alpha)\\[1.5ex]
&&
\;-\;
L^{s-1} \circ \Lambda \,(L \alpha)
\;+\;
(m-p)\,L^{s-1}\alpha
\\[1.5ex]
&=&
((s-1)(m-p-s) \,+\, (m-p))\,L^{s-1}\alpha
\\[1.5ex]
&=&
s\,(m-p-s+1)\,L^{s-1}\alpha \ . 
\eea
The formula for $\,\Lambda^r L^s\,$ in the case $r>1$ then follows 
by applying the commutator relation several times. \qed

On a K{\"a}hler manifold, any $p$--form $\,u$ with $\,p\le m$ has a unique
decomposition, $\,u = u_0 + L u_1 + \ldots + L^l u_l$, where
the $\,u_i$'s are primitive forms. This is usually called the {\it Weyl 
decomposition} of $\,u$. From now on we suppose that $p\le
m$, otherwise we just replace $u$ by its Hodge dual $*u$. Applying
$\,\Lambda L\,$ 
and using the Lemma~\ref{elementary} in the case $\,r=1\,$ we obtain
$$
\Lambda\,L(u) 
\;=\;
(m-p)\,u_0
\;+\;
2\,(m-p+1)\L u_1
\;+\;\ldots\;+\;
(l+1)\,(m-p+l)\,L^lu_l \ .
$$
It is easy to see that coefficients in the above sum are all different.
Indeed, $(s+1)(m-p+s)=(r+1)(m-p+r)$ if and only if
$(s-r)(m-p+1+s+r)=0$, recall that we assumed that $p\le m$. Hence, a
$p$--form $\,u\,$ with $\,p\le m\,$ is 
an eigenvector of $\,\Lambda L\,$ if and only $\,u=L^i\alpha$, for
some primitive form $\,\alpha$ and in that case the corresponding
eigenvalue is $(i+1)(m-p+i)$. 
 
Lemma \ref{myh} thus implies that $du=L^sv$ and $\d u=L^rw$ for some
primitive 
forms $v$ and $w$, and moreover
$\,(2m-p-2)(p+2)=4(s+1)(m-p+s-1)\,$ and
$\,p(2m-p)=4(r+1)(m-p+r+1)$, which have the unique solutions
$\,2s=p\,$ and $\,2r=p-2$. Hence, denoting $p=2k$, we have
$$
du \;=\; L^k v
\qquad\qquad
\d u\;=\; L^{k-1}w
$$
for some vectors $v$ and $w$. Moreover, we see that the
Weyl decomposition of $\,u\,$ has the form
$\,u=u_0 + \ldots + L^ku_k$, where $\,u_k$ has to be a function.

\begin{elem}\label{jh}
Let $\,u\,$ be a $2k$--form with $Ju$ parallel and such that
$\,du = L^kv\,$ and 
$\,\d u = L^{k-1}w$, for some vectors $v$ and $w$. Then
$$
w \;=\; \tfrac{k(2m-2k+1)}{2k+1}\; Jv \.
$$
\end{elem}  
\proof
Applying $J$ to the equation $\,du = L^kv\,$ leads to
$\,d^cu = Jdu = L^k(Jv)$. Thus, contracting with the
K{\"a}hler form yields
$$
\Lambda(d^cu) 
\;=\;
 d^c\Lambda u \;+\; \d u
\;=\;
\Lambda \, L^k (Jv)
\;=\;
k\,(m-k)\,L^{k-1}(Jv) \ .
$$
From Proposition~\ref{differentials} we have the equation 
$\,\d u = - \mu_1 d^c \Lambda u$. Hence,
$$
\d u
\;=\;
(1 \,-\, \tfrac{1}{\mu_1})^{-1}k(m-k) L^{k-1}(Jv) \ ,
$$
which proves the lemma after substituting the value of
$\mu_1$. \qed

In the next step we will show that in the Weyl decomposition
of $\,u\,$ only the last two terms, i.e. $\,L^{k-1}u_{k-1}\,$
and $\,L^ku_k$, may be non--parallel. Indeed we have
\bea
\nabla_X u
&=&
\nabla_X u_0 \;+\; L(\nabla_Xu_1)\;+\;\ldots\;+\;L^k(\nabla_Xu_k)
\\[1.5ex]
&=&
\tfrac{1}{p+1}\,X \es L^kv \;-\; \tfrac{1}{n-p+1}\,X \we L^{k-1}w
\\[1.5ex]
&=&
L^{k-1}\, ( \tfrac{k}{p+1}\,JX \we v \,-\, \tfrac{1}{n-p+1}\,X \we w)
\;+\;
L^k(\tfrac{1}{p+1}\,\la x, \, v \ra)
\\[1.5ex]
&=&
L^{k-1} \, \omega_1 \;+\; L^k \omega_2  \ ,
\eea
for some primitive forms $\,\omega_1,\,\omega_2$. Hence,
comparing the first line with the last, implies that the
components $\,u_0,\,u_1,\ldots,u_{k-2}\,$ have to be parallel.
But then we can without loss of generality assume that they are
zero, i.e. for a twistor $p$--form $\,u$, with $\,Ju\,$ parallel,
we have
$$
u \;=\; L^{k-1}u_{k-1}\;+\; L^k u_k
$$
where $\,u_{k-1}\,$ is a primitive 2--form and $\,u_k\,$ is a
function,~i.e. $\,p=2k$.

Our next aim is to translate the twistor equation for $u$ into 
equations for $u_{k-1}$ resp. $u_k$. This leads naturally to 
the following definition, which we will need for the rest of
this paper.
\begin{ede} \label{spec} A {\em special} 2--form on a K{\"a}hler manifold $M$ is a
  primitive 2--form $\f$ of type (1,1) satisfying the equation
\begin{equation}\label{invariant1}
\nabla_X \, \f
\;=\;
\gamma \we JX \;-\; J\gamma \we X
\;-\;
\tfrac{2}{m}\,\gamma(X)\,\omega \ ,
\end{equation}
for some 1--form $\gamma$, which then necessarily equals 
$\,\tfrac{m}{2(m^2-1)}\d^c \f$.
\end{ede}

Recall that a 2--form $u$ is of type (1,1) if and only if it is
$J$ invariant,~i.e. $Ju=0$. 
As a first property of such special forms we obtain

\begin{elem} If $\f$ is a special 2--form on a compact K{\"a}hler
  manifold $M$ of dimension $m>2$, then $\dec \f$ is exact.
\end{elem}
\begin{proof}
Taking the wedge product with $X$ in (\ref{invariant1}) and summing
over an orthonormal basis yields
$d\f=2\tfrac{m-1}{m}\gamma=\tfrac{1}{m+1}\d^c \f$. 
It follows that $\,Ld\d^c\f = 0\,$ and, since
$\,L\,$ is injective on 2--forms, we conclude that $\,\d^c\f\,$ is
closed. Hence, 
$\,\d^c\f = h + df\,$ for some function $\,f\,$ and a harmonic 1--form
$\,h$. We have to show that $\,h\,$ vanishes. First, note that 
$\,h\,$ is in the kernel of $\,d^c\,$ and $\,\d^c\,$ since the manifold is
K{\"a}hler. Computing its $\,L^2$-norm we obtain
\bea
(h,\,h)
&=&
(h, \, h + df)
\;=\;
(h,\, \d^c\f)
\;=\;
(d^c h,\, \f)
\;=\;
0 \ .
\eea
\r

\begin{ede} The function $f$ given by the lemma above will be called
the {\em generalized trace} of the special form $\f$. It is only 
defined up to a constant. 
\end{ede}

We can now state the main result of this section

\begin{ath}\label{main} Let $u$ be a form of degree $p$ on a compact
  K{\"a}hler manifold $M^{2m}$ and suppose that $n-2\ge p\ge 2$ and $p\ne
  m$. Then $u$ is a 
  twistor form if and only if there exists a special form $\f$ whose
  generalized trace is $f$ and a positive integer $k$ such that $p=2k$ and 
$$
u\;=\;L^{k-1}\f\;-\;\tfrac{m-p}{p(m^2-1)}L^kf\;+\;\hbox{parallel form}.
$$
The same statement is valid for $p=m$ under the additional assumption
that $Ju$ is parallel.
\end{ath}
\begin{proof}
Consider first the case $p\le m$.
Using the notations from Lemma \ref{jh}, the twistor equation for $u$
reads
\begin{equation}\label{ui}
\n_Xu_{k-1}\;+\;L\n_Xu_k
\;=\;
\tfrac{k}{p+1}JX \we v \;-\; \tfrac{1}{n-p+1}X \we w .
\end{equation}
An equality between 2--forms is equivalent with the equality of
their primitive parts and that of their traces with respect to the
K{\"a}hler form.
The equality of the primitive parts in (\ref{ui}) is equivalent to
$$
\n_Xu_{k-1}
\;=\;
(\tfrac{k}{p+1}JX \we v \,-\, \tfrac{1}{n-p+1}X \we w)
\;-\;
\tfrac{1}{m}\,
\Lambda \,(\tfrac{k}{p+1}JX \we v \,-\, \tfrac{1}{n-p+1}X \we w)
\, \omega \ ,
$$
where $\,\omega \,$ is the K{\"a}hler form. Now,
$\,\Lambda(JX \we v) = - \la v, X \ra  \,$
and
$\, \Lambda (X \we w ) - \la Jw, X \ra \,$ so we obtain
\bea
\firstline{.1cm}{
\n_Xu_{k-1}
\;=\;
\tfrac{k}{p+1}\, JX \we v \,-\, \tfrac{1}{n-p+1} \, X \we w
\;+\;
\tfrac{1}{m}\,
(\tfrac{k}{p+1}\,\la v, X \ra
\,-\,
\tfrac{1}{n-p+1}\, \la Jw, X \ra   )\, \omega}
\\[1.5ex]
&=&
\tfrac{k}{p+1}\, JX \we v \,-\, \tfrac{1}{n-p+1}
\,\tfrac{k(n-p+1)}{p+1}\, X \we Jv
\;+\;
(\tfrac{k}{m(p+1)}
\,-\,
\tfrac{1}{m(n-p+1)} \tfrac{k(n-p+1)}{p+1}) \,  \la v, X \ra \,\omega
\\[1.5ex]
&=&
\tfrac{k}{p+1}\, JX \we v \,-\, \tfrac{k}{p+1} \, X \we Jv
\;+\;
\tfrac{2k}{m(p+1)}\,\la v, X \ra \, \omega
\\[1.5ex]
&=&
\gamma \we JX \;-\; J\gamma \we X \;-\;
\tfrac{2}{m}\,\gamma(X)\,\omega \ ,
\eea
where $\,\gamma\,$ is the 1--form defined by 
$\,\gamma(X) = -\tfrac{k}{p+1} \la v, X \ra$. Contracting with $JX$
shows that $\gamma=\tfrac{m}{2(m^2-1)}\d^c u_{k-1}$, so
$v=\tfrac{(p+1)m}{p(m^2-1)}\d^c u_{k-1}$.

The second part of (\ref{ui}) consists in the equality of the traces,
which is equivalent to
$$
X(u_k)
\;=\;
-\,\tfrac{p}{m(p+1)}\,\la v, X \ra
\;+\;
\tfrac{1}{p+1}\, \la v, X \ra
\;=\;
\tfrac{m-p}{m(p+1)}\,X \es v \ ,
$$
i.e.
$$
d u_k \;=\; \tfrac{m-p}{m(p+1)}\, v \ =-\tfrac{m-p}{p(m^2-1)}\d^c u_{k-1}.
$$
We have shown that $u$ is a twistor form if and only if $u_{k-1}$ is
a special form and $u_k$ is $-\tfrac{m-p}{p(m^2-1)}$ times its
generalized trace.

A priori, this only proves the theorem for $p\le m$, but because of
the invariance of the hypothesis and conclusion of the theorem with
respect to the Hodge duality, the theorem is proved in full generality.
\r

Specializing the above result for $k=1$ yields the following
characterization of twistor 2--forms also obtained in
\cite{acg}. Note that the result that we obtain is more
complete than that of  \cite{acg} since we do not assume the
invariance of the 2--form $u$.

\begin{elem}\label{character}
Let $(M^{2m},\,g,\,J,\,\omega)$ be a K{\"a}hler manifold. Then a
2--form $\,u\,$ is a twistor form if and only if there
exists a 1--form 
$\gamma$  with
\begin{equation}\label{invariant}
\nabla_X \, u
\;=\;
\gamma \we JX\;-\; J\gamma \we X
\;-\;
\gamma(X)\,\omega \ .
\end{equation}
The 1--form $\gamma$ is necessarily equal to
$\gamma = \tfrac{1}{2m-1}\,J\d u$. Moreover, 
$\,\gamma = -\,\tfrac{1}{m-2}\,d\, \la u, \omega \ra \,$ provided that $\,m>2$.
\end{elem}
\begin{proof}
First of all we can rewrite equation~(\ref{invariant}) as
$\,\nabla_X u = - X \es(\gamma\we \omega) + X \we J\gamma $.
Contracting (resp. taking the wedge product) in (\ref{invariant}) with 
$X$ and summing over an orthonormal frame 
$\,\{e_i\}\,$ we obtain
$$
J\gamma \;=\; -\,\tfrac{1}{2m-1}\,\d u
\qquad \mbox{and} \qquad
\gamma \we \omega \;=\; -\,\tfrac{1}{3}\,du \ .
$$
Substituting this back into~(\ref{invariant}) yields the defining
equation for a twistor 2--form, i.e.:
$$
\nabla_X\,u
\;=\;
\tfrac{1}{3}\,X \es du \;-\;
\tfrac{1}{2m-1}\,X \we \d u \ .
$$
Conversely, Theorem \ref{main} shows that $u=u_0+Lf$, where 
$$\nabla_X \, u_0
\;=\;
(\gamma \we JX \;-\; J\gamma \we X) 
\;-\;
\tfrac{2}{m}\,\gamma(X)\,\omega \ ,
$$
and $df=-\tfrac{m-2}{2(m^2-1)}\dec u_0=-\tfrac{(m-2)}{m}\gamma$, so
clearly $u$ satisfies (\ref{invariant}).
\r

This characterization of twistor 2--forms in particular implies
that for $\,m>2\,$ special forms are just the primitive 
parts of twistor 2--forms and vice versa. 
In the remaining part of this section we will describe a similar
relation between twistor forms and Hamiltonian 2--forms. Using the 
results of \cite{acg} we can thus produce many examples of
non--parallel twistor forms.

\begin{ede}
An invariant 2--form $\,\psi\,$ is called {\it Hamiltonian} if there
is a function $\,\sigma\,$ such that
$$
\nabla_X \psi
\;=\;
\tfrac{1}{2}\,
(d\sigma \we JX  \;-\; Jd\sigma \we X  )
$$
for any vector field $X$. When $m=2$ one has to require in addition
that $\,J\grad(f)\,$ is a Killing vector field.
\end{ede}

It follows immediately from the definition that 
$\,d\sigma = d \la \psi, \omega \ra$. Hence, one could
without loss of generality replace $\,\sigma\,$
with $\,\la \psi, \omega \ra$.

Any special 2--form $\f$ with generalized trace $f$ defines
an affine line $\f + \RM f\omega$ in the space of 2--forms modulo constant
multiples of the K{\"a}hler form. 
Lemma~\ref{character} shows that this line contains 
a twistor 2--form and it is not difficult to see 
that it contains a unique closed form and also an unique coclosed
form. Indeed, for some real number 
$x$, $d(\f+xf\o)=0$ is equivalent to $-\frac{1}{m+1}L\dec\f+xLdf=0$,
i.e. $x=\frac{1}{m+1}$, and $\d(\f+xf\o)=0$ is equivalent to
$\d\f=x\dc f=xJdf=xJ\dec\f=-x\d\f$, i.e. $x=-1$. 
The following proposition
shows that this affine line also contains a Hamiltonian 2--form.

\begin{epr}
Let $(M^{2m},\,g,\,J,\,\omega)$ be a K{\"a}hler manifold with a
Hamiltonian 2--form $\,\psi$. Then 
$\,u:= \psi - \tfrac{\la \psi, \omega \ra}{2} \,\omega\,$ 
is a twistor 2--form. Conversely,
if $\,u\,$ is a twistor 2--form and $\,m>2$, then 
$\,\psi :=u-\tfrac{\la u, \omega \ra }{m-2}\,\omega\,$ is
a Hamiltonian 2--form.
\end{epr}
\proof
Let $\,\psi\,$ be a Hamiltonian 2--form and let $\,u\,$ be defined as
$\,u:= \psi - \tfrac{\la \psi, \omega \ra}{2} \,\omega$, then
\bea
\nabla_X\,u
&=&
\nabla_X\,\psi 
\;-\;
\tfrac{1}{2}\,d \la \psi, \omega \ra (X)\,\omega\\[1.5ex]
&=&
\tfrac{1}{2}\,
\left(
d\sigma \we JX \;-\; J d\sigma \we X
\right) 
\;-\;
\tfrac{1}{2}\,d\sigma (X)\,\omega \ .
\eea
Thus, the invariant 2--form $\,u\,$ satisfies the equation~(\ref{invariant})
of Lemma~\ref{character} with, $\,\gamma:=\tfrac{1}{2}\,d\sigma$,
and it follows that $\,u\,$ is a twistor 2--form. Conversely, starting
from a twistor 2--form $\,u\,$ and defining 
$\,\psi :=u-\tfrac{\la u, \omega \ra }{m-2}\,\omega$, we can use the
characterization of Lemma~\ref{character} to obtain
\bea
\nabla_X \psi 
&=& 
\gamma \we JX \;-\; J\gamma \we X \ .
\eea
Since for $\,m>2\,$ we have 
$\,\gamma =- \tfrac{1}{m-2}d\la u, \omega \ra 
= \tfrac{1}{2}d\la \psi, \omega \ra$, which implies that $\,\psi\,$
has to be a Hamiltonian 2--form.
\qed

The situation is somewhat different in dimension 4. If $\,\psi\,$
is a Hamiltonian 2--form then 
$u = \psi -\tfrac{\la \psi, \omega\ra }{2}\omega =u_0\,$
is an anti--self--dual twistor 2--form (or equivalently, invariant and
primitive). Conversely, if we start 
with an anti--self--dual twistor 2--form $\,u_0\,$ and ask for which
functions $\,f\,$ the invariant 2--form $\,\psi:=u_0 + f \omega\,$ is
a Hamiltonian 2--form we have
\begin{elem}\label{nochda}
Let $\,(M^4,\,g,\,J)\,$ be a K{\"a}hler manifold with an anti-self-dual
twistor 2--form $\,u_0$. Then $\,\psi:=u_0 + f \omega\,$ is Hamiltonian 2--form if
and only if
\begin{equation}\label{hamcond}
\d u_0 
\;=\;
-\,3\,Jdf 
\end{equation}
and $\,\d u_0\,$ is dual to a Killing vector field.
In particular, this is the case on simply connected Einstein manifolds.
\end{elem}
\proof
Since $f=\tfrac{\la \psi, \omega \ra}{2}$ we conclude from the
definition 
that $\,\psi\,$ is
a Hamiltonian 2--form if and only if 
$\nabla_X \psi = df \wedge JX - Jdf \wedge X$
for any vector field $X$. On the other hand $u_0$ is assumed to be an invariant 
twistor 2--form. Hence, from Lemma~\ref{character}
$u_0$ satisfies
$ \nabla_X\,u_0 = - X \es (\gamma \wedge \omega)- J\gamma \wedge X$,
where $\gamma=\tfrac{1}{3}J\d u_0$.
This implies that $\psi = u_0 + f \omega$ is a Hamiltonian 2--form if
and only if
$$
- X \es (\gamma \wedge \omega)- J\gamma \wedge X
\;=\;
df \wedge JX - Jdf \wedge X - df(X)\omega
\;=\;
- X \es (df \wedge \omega)- Jdf\wedge X .
$$
This is the case if and only if $\,\gamma = df\,$, or equivalently if
$\,\d u_0 \;=\;-\,3\,Jdf$. If the complex dimension is 2, the
definition of Hamiltonian 2--forms made the additional requirement
that $\,Jd\sigma$, which here is proportional to $\,\d u_0$, is dual
to Killing  vector field.
Hence, it remains to show that on simply connected Einstein manifolds the 
equation~(\ref{hamcond}) has a solution such that $\,\d u_0\,$ is dual 
to a Killing vector field. Let $X$ be any Killing vector field on an arbitrary
K{\"a}hler manifold, then 
$$
0 
\;=\;
L_X\,\omega
\;=\;
d\,X\es \omega \;+\; X\es d\,\omega
\;=\;
d\,J\,X \ ,
$$
i.e. the 1--form $\,JX\,$ is closed. Hence, since the manifold is
simply connected, there exists some function $\,f\,$ with
$\,JX = df\,$ and also $\,X = -Jdf$. Finally, it is easy to see  
that on Einstein manifolds, for any twistor 2--form $\,u_0$, 
the 2--form $\,\d u_0\,$ is dual to a Killing vector field (cf. \cite{uwe}).
 \qed

Hamiltonian 2--forms have the following remarkable property
(c.f.~\cite{acg}). 
If $\,\psi\,$ is Hamiltonian, and if $\,\sigma_1,\ldots,\sigma_m\,$
are the elementary symmetric functions in the eigenvalues of
$\,\psi\,$ with respect to the K{\"a}hler form $\,\omega$,
then all vector fields $\,K_j = J\grad(\sigma_j)\,$ are Killing.
Furthermore, the Poisson brackets $\,\{\sigma_i, \sigma_j\}\,$
vanish, which implies that the Killing vector fields
$\,K_1,\ldots,K_m\,$ commute. If the Killing vector fields
are linearly independent, then the K{\"a}hler metric is toric.
But even if they are not linearly independent one has further
interesting properties, which eventually lead to a  complete local 
classification of K{\"a}hler manifolds with Hamiltonian 2--forms
in \cite{acg}.
The most important sources of K{\"a}hler manifolds with
Hamiltonian 2--forms are weakly Bochner--flat K{\"a}hler manifolds
and (in dimension greater than four) K{\"a}hler manifolds
which are conformally-Einstein.
A K{\"a}hler manifold is weakly Bochner--flat if its Bochner tensor
is coclosed, which is the case if and only if the normalized 
Ricci form is a Hamiltonian 2--forms.
The examples include some Hirzebruch surfaces and the complex 
projective spaces.


\section{The middle dimension}

In this section we study twistor forms of degree $m$ on compact
K{\"a}hler manifolds of real dimension $n=2m$. This case is very special
thanks to the following 

\begin{elem}\label{parallel} 
Let $\L^mM=L_1\oplus\ldots\oplus L_r$ be a decomposition
  of the bundle of $m$--forms in parallel subbundles and
  $u=u_1+\ldots+u_r$ be the corresponding decomposition of an
  arbitrary $m$--form $u$. Then $u$ is a twistor form if and only if
  $u_i$ is a twistor form for every $i$.
\end{elem}
\proof
In the case of $m$--forms on  $2m$-dimensional manifolds 
Proposition~\ref{integrabl} provides a characterization of
twistor forms similar to that of Killing forms. A $m$--form
$u$ is a twistor form if and only if
\begin{equation}\label{midchar}
\Delta \, u \;=\; \frac{m+1}{m}\,q(R)\,u
\end{equation}
Given a decomposition of the form bundle $\Lambda^mM$ into 
parallel subbundles we know that the Laplace operator $\Delta$
as well as the symmetric endomorphism $q(R)$ preserve this
decomposition. Hence, equation~(\ref{midchar}) can be projected
onto the summands,~i.e. (\ref{midchar}) is satisfied for 
$u=u_1+\ldots +u_r$ if and only it is satisfied for all summands
$u_i$, which proves the lemma. \qed

In Section~\ref{kaehler} we defined special 2--forms, which turned
out to be the main building block for twistor $p$--forms with $p\neq m$.
In dealing with twistor $m$--forms on 2$m$--dimensional K{\"a}hler manifolds
we have to introduce the following

\begin{ede}
A special $m$--form on a 2$m$--dimensional K{\"a}hler manifold is a 
$m$--form $\psi$ of type $(1,m-1)+(m-1,1)$ satisfying for all vector fields
$X$ the equation
\begin{equation}\label{spec-m}
  \nabla_X\, \psi 
  \;\;=\;\;
  (m-1)\,\tau \wedge JX \;-\; \,J\tau \wedge X 
  \;-\; (m-1)\,(X\es \tau)\wedge\omega \ , 
\end{equation}
for some 1--form $\tau$, which then necessarily equals 
$\,\frac{1}{m^2-1}\d^c \psi$.
\end{ede}

Note that for $m=2$, we retrieve the definition of special 2--forms on
4--dimensional manifolds (Definition \ref{spec}).

\begin{epr} \label{middle}
  Let $u$ be a $m$--form on a compact K{\"a}hler manifold
  $M^{2m}$. Then $u$ is a twistor form if and only if $m=2k$ and there 
  exists a special 2--form $\f$ and a special $m$--form $\psi$ such that
$$
u\;\;=\;\;L^{k-1}\,\f\;+\;\psi\;+\;\hbox{parallel form}.
$$
\end{epr}
\proof
Lemma \ref{parallel} shows that we can assume  $u=L^ku_k$ for some 
primitive form $u_k$. Suppose
$k\ge 1$, then it follows from Proposition~\ref{jlu} that 
$\,J \Lambda u = J \Lambda L^k u_k = 0$. Thus Lemma~\ref{elementary}
implies $\,JL^{k-1}u_{k}=0$. Since $\,L^{k-1}\,$ is injective on
$\,\Lambda^{m-2k}$, we get $\,Ju_k=0\,$ and eventually $\,Ju=0$.
For $J$--invariant twistor forms the conclusion of
Theorem~\ref{main} also holds in the case  $p=m$. Hence, $\, u= L^{k-1}\phi + 
parallel \,form$, where $\phi$ is a special 2--form.

It remains to treat the case $\,k=0$,~i.e. the case where $\,u\,$ is a
primitive $m$--form. Contracting the twistor equation with the K{\"a}hler
form yields
\begin{equation}\label{type1}
X \es \d^c u \;+\; JX \es \d u \;\;=\;\; 0
\end{equation}
for all vector fields $X$. After wedging with $X$ and $JX$ and 
summing over an orthonormal basis, this gives
\begin{equation}\label{type}
J\,\d\,u \;=\; (m-1)\,\d^c\,u
\qquad\mbox{and}\qquad
J\,\d^cu \;=\;-\,(m-1)\,\d\,u \ .
\end{equation}
Now, $\,\ast u\,$ is again a twistor $m$--form and the above argument 
shows that we can suppose $\,\ast u\,$ to be primitive. (Otherwise,
$\,\ast u = L^{k-1}\phi + parallel \,form$, for a special 2--form
$\,\phi$, which in particular is self-dual,~i.e. $\,u = L^{k-1}\phi
+ parallel\, form$). Applying the Hodge star operator to $\,\Lambda(\ast
u)=0\,$ immediately implies $Lu=0$. This shows that $\,du=-L\d^cu$,
so the twistor equation becomes
\begin{eqnarray}
\nabla_X \, u&=& 
\frac{1}{m+1}\,(\, 
-\, X \es L\,\d^c\,u \;-\; X \we \d \,u\,) 
\label{twistor}\\[1.5ex]
&=&
\frac{1}{m+1}\,(\,
-\,JX \we \d^cu 
\;-\;L\,(X\es \d^cu)
\;+\;
\frac{1}{m-1}\,X \we J\d^cu 
\,)
\ ,\nonumber
\end{eqnarray}
which is just the defining equation~(\ref{spec-m}) of a special 
$m$--form with $\tau = \,\frac{1}{m^2-1}\d^cu$. According to our
definition of special $m$--forms, we still have to show that (up to
parallel forms)
$u$ is of type $(m-1,1)+(1,m-1)$. Equivalently we will show this
for $\nabla_X u$, where $X$ is any vector field. Recall that $u$
is of type $(m-1,1)+(1,m-1)$ if and only if $\,J^2u=-(m-2)^2u$. We
will compute $\,J^2(\nabla_X u)\,$ using the twistor equation. 
Equation~(\ref{type}) implies that $\,\d u\,$ is of type 
$(m-1,0)+(0,m-1)$,~i.e $\,J^2(\d u)=- (m-1)^2\d u$. Moreover,
using (\ref{type1}) and (\ref{type}) we obtain
\bea
J^2(X \we \d u)
&=&
(-1-(m-1)^2)\,X \we \d u \;+\; 2\,JX \we J\d u\qquad
\\[1.5ex]
&=&
(-1-(m-1)^2)\,X \we \d u \;+\;2\,(m-1)\,JX \we \d^c u
\eea
and similarly 
\bea
\firstline{1.5cm}{
J^2(X\es L\d^c u)
\;=\;
(-1-(m-1)^2)\,X  \es L\d^cu \;+\; 2\,JX \es JL\d^c u}
\\[1.5ex]
&=&
(-1-(m-1)^2)\,X  \es L\d^cu \;-\;2\,(m-1)\,JX \es L\d u
\\[1.5ex]
&=&
-(m-2)^2\,X \es L\d^cu \;+\; 2\,(m-1) X \we \d u
\;-\; 2(m-1)\, JX \we \d^c u \ .
\eea
Adding these two equations implies
$\,J^2(\nabla_X u)=-(m-2)^2\nabla_Xu\,$
after using the twistor equation as written in~(\ref{twistor}). \qed
 
As an application of Proposition~\ref{middle} we see that any twistor 
2--form on a 4-dimensional K{\"a}hler manifold has to be invariant and
primitive,~i.e. any twistor 2--form has to be a special 2--form. In this
case we know from Lemma~\ref{nochda} that any Hamiltonian 2--form 
gives rise to a twistor 2--form and vice versa, under some additional
conditions. This clarifies the situation of twistor 2--forms in
dimension~4 and shows in particular that one can exhibit many examples. 
However,  one can  show that any special 
$m$--form ($m\ge3$) on a K{\"a}hler-Einstein manifold has to be
parallel. Nevertheless for the moment it  
remains unclear whether one can exclude these forms in general.


\section{Twistor forms on the complex projective space}

In this section we describe the the construction of twistor forms on 
the complex projective space (c.f.~\cite{acg}).
Let $\,M = \CM P^m$ be equipped with the Fubini-Study metric and the 
corresponding K{\"a}hler form $\,\omega$. Then the Riemannian curvature 
is given as
$$
R_{X,\,Y}Z
\;=\;
-\,(X \wedge Y \;+\; JX \wedge JY)\,Z
\;-\;
2\,\omega (X,\,Y)\,JZ
$$
for any vector fields $\,X, Y, Z$. This implies for the Ricci 
curvature $\,\Ric = 2(m+1) \id $. Let $\,K\,$ be any Killing vector 
field on $\,\CM P^m$. Then there exists a function $f$ with
$\,\Delta f = 4(m+1)f\,$ and $\,K = J \grad(f)$, i.e. $\,f\,$ is an
eigenfunction of the Laplace operator for the first non--zero
eigenvalue. Now, consider the 2--form 
$\,\phi:= dK=dJdf = dd^c(f)$. Since $\,K\,$ is a Killing vector field
it follows:
\bea
\nabla_X \,\phi
&=&
\nabla_X(dK) \;=\;2\,\nabla_X(\nabla K)
\;=\;
2\,\nabla^2_{X,\,\cdot}\,K
\;=\;
-\,2\,R(K,\,X)
\\[1.5ex]
&=&
-\,2\,
(df \we JX \;-\; Jdf \we X)
\;-\;
4\,df(X)\,\omega \ .
\eea
It is clear that $\,\phi\,$ is an invariant 2--form and an eigenform of the
Laplace operator for the minimal eigenvalue $\,4(m+1)$. A small modification
of $\,\phi\,$ yields a twistor 2--form. Indeed, defining
$\,{\hat \phi}:= \phi + 6 f\omega\,$ one obtains
$$
\nabla_X \,{\hat \phi}
\;=\;
-\,2\,
(df \we JX \;-\; Jdf \we X)
\;+\;
2\,df(X)\,\omega \ .
$$
Using Lemma~\ref{character} for $\,\gamma:=-2\,df\,$ one concludes that
$\hat \phi$ is a twistor 2--form. It is not difficult to show
that indeed any invariant twistor 2--form on $\,\CM P^m\,$ has 
to be of this form. Summarizing the construction  one has

\begin{epr}[\cite{acg}]
Let $K = J \grad(f)$ be any Killing vector field on the complex projective
space $\CM P^m$ then
$$
{\hat \phi}
\;:= \;
d d^c (f)\;+\; 6\,f\,\omega
\;=\;
(d d^c (f))_0\;+\; \tfrac{2m-4}{m}\,f\,\omega
$$
defines an invariant, non--parallel twistor 2--form. Moreover, 
in dimension four, $\,{\hat \phi}\,$ is a primitive
2--form,~i.e. $\,{\hat \phi=(d d^c (f))_0}$.
\end{epr}


\section{Concluding Remarks}

For the convenience of the reader, we resume here the main results of
the paper, as well as the remaining open questions. 

If $u$ is a twistor $p$--form on a compact K{\"a}hler manifold $M^{2m}$
with $2\le p \le n-2 $ and $p\ne m$, then it turns out that $Ju$ has to be parallel (Corollary
\ref{jin}). In this case $u$ is completely determined (modulo parallel forms)
by a special 2--form and its generalized trace (Theorem \ref{main}).
In particular, the degree of $u$ has to be even. 

If $p=m$, the $J$--invariant part of $u$ can still be characterized
(up to parallel forms) by special 2--forms, but the complete
classification of non--parallel twistor forms can only be obtained up to
the hypothetical existence of {\it special $m$--forms} of type
$(m-1,1)+(1,m-1)$. At this moment we have no examples of such forms
for $m>2$, and moreover we were able to show that they don't exist on compact
K{\"a}hler--Einstein manifolds. Note that for $m=2$, the definition
of special $m$--forms coincides with the previous one for special
2--forms on 4--manifolds.

For $m>2$, each special 2--form $\f$ with generalized trace $f$
(defined by
$df=\dec\f$) determines an affine line $\{\f_x=\f+xf\o\}$ in the space of
two--forms, which contains distinguished elements~: $\f_x$ is
Hamiltonian for $x=\frac{1}{m^2-1}$, closed for $x=\frac{1}{m+1}$, coclosed for
$x=-1$ and a twistor form for $x=-\frac{m-2}{2(m^2-1)}$. 

For $m=2$, the picture is slightly different since the generalized
trace of a special 2--form (and thus the affine line above) are no
longer defined. The results are as follows~: every twistor 2--form is
(up to parallel forms) primitive and invariant
(i.e. anti--self--dual), and defines a special 2--form. This special
2--form induces a Hamiltonian 2--form if and only if its
codifferential is a Killing vector field. Conversely, a Hamiltonian
2--form always defines the affine line above, and in particular its
primitive part is simultaneously a special 2--form and a twistor
form. There are examples of twistor 2--forms on 4--dimensional
manifolds which do not come from Hamiltonian 2--forms \cite{apc}. The
complete local classification of Hamiltonian 2--forms was obtained in 
\cite{acg1} for $m=2$ and \cite{acg} in the general case. The same
references contain examples of compact K{\"a}hler manifolds with
non--parallel Hamiltonian 2--forms.


 \labelsep .5cm

\end{document}